\documentclass{article}

\usepackage{amsmath,amssymb}
\usepackage{color}
\usepackage{graphicx}

\title{STUDENTS AS PARTNERS IN CURRICULAR DESIGN: \\ CREATION OF STUDENT-GENERATED CALCULUS PROJECTS}

\author{Steve Cohen\\
Department of Mathematics and Actuarial Science\\
Roosevelt University\\
Chicago, IL 60605, USA\\
scohen@roosevelt.edu
\and
B\'arbara Gonz\'alez-Ar\'evalo \\
Department of Mathematics\\
Hofstra University\\
Hempstead, NY 11549, USA\\
barbara.p.gonzalez@hofstra.edu
\and
Melanie Pivarski \\
Department of Mathematics and Actuarial Science\\
Roosevelt University\\
Chicago, IL 60605, USA\\
mpivarski@roosevelt.edu
}

\providecommand{\keywords}[1]{\textbf{\textit{Keywords:}} #1}

\newcommand{\AmSLaTeX}{$\cal A$\kern-.1667em\lower.5ex\hbox{$\cal
M$}\kern-.125em $\cal S$-\LaTeX}

\newcommand{\RU}{Roosevelt University}

\begin{document}

\maketitle

\begin{abstract}
For the past five years advanced students have developed projects
for our redesigned Calculus II classes to use. Our student project designers are
often mathematically early in their careers, and so this provides them with an
opportunity to create and explore new mathematics while giving us the ability
to involve students of all levels in research projects. It also gives them the chance to
present their work at conferences and local talks.

Our course redesign affected three groups of students: ones taking the
class, ones designing projects for the course, and embedded tutors. This qualitative study examines how the second and third groups of students benefited from their
experiences and how we can modify our program to improve it. Evidence was gathered from interviews, surveys, and observation of student research work and its implementation in the classroom.  We found that tutors reported more confidence in their knowledge of calculus and insights into teaching it, and project designers experienced similar benefits to that of a traditional research experience.

\end{abstract}


\section{INTRODUCTION}

The extensive use of undergraduate research in mathematics is a fairly recent one, dating back to the 1980s with the widespread introduction of the NSF-funded Research Experience for Undergraduates programs \cite{MAA}.  Most of these experiences are designed for advanced undergraduates who are in their junior or senior years, and they are often used to help prepare these students for graduate study.  By using undergraduate students to develop projects for use in a Calculus II classroom, we also allow students who are freshmen and sophomores the chance to work on research.  The use of this research is clear; our students are motivated by helping their peers learn.  They also enjoy the flexibility involved in choosing a topic that allows them to explore their own interests.

The use of undergraduates as embedded peer tutors is common; see e.g.\cite{BenefitsTA, PeerTutor}.  Tutors attend most classes, and depending on the instructor, work with students during the class sessions.  This gives students a chance to connect more deeply with the material, increasing their calculus skills as well as their empathy for their fellow students.  

In order to avoid ambiguity, we use ``embedded tutors'' or simply ``tutors'' to refer to the embedded peer tutors and ``project designers'' or ``student researchers'' to refer to students who, after completing Calculus II themselves, worked on researching a project for use in a future Calculus II course.  We refer to students who were currently taking the course simply as calculus students.

In section \ref{connecting}, we describe our Calculus II course and the role of the tutors and project designers.
In section \ref{curricular}, we discuss definitions of student research that occur in the literature and how these connect to our curricular design.
In section \ref{methodology} we describe the methodology used in our study.
In section \ref{results}, we delve into the results of the study, providing and elaborating on themes found in the student responses.

\section{CONNECTING STUDENTS TO THE COURSE}\label{connecting}
At \RU, semester-long projects became a regular part of all sections of Calculus II in spring of 2010 \cite{Roosevelt}.  Calculus students use
projects as a way to explore STEM majors that use calculus, to develop library research skills, and to develop communication skills.  Since
then these classes have had embedded undergraduate tutors who have
had Calculus II at \RU or another institution. Tutors attend class at
least once per week and help students in and out of class. The Calculus
II classes each have between 9 and 30 students, and there is one peer
tutor assigned to each section. Peer tutors schedule weekly office hours,
assist with computer activities, help students with homework, and help
groups with their projects. The level of actual involvement varies based on the instructor and the particular students in the class.

Starting in summer 2011 undergraduate research students were given
the opportunity to work on designing materials to use for class projects.
Their work involves picking a topic of civic importance, finding appropriate data sources, considering issues related to calculus, and linking these
together. There are many possible outcomes for these projects: use in
a Calculus II class, honors theses, research talks, and as starter ideas
for more advanced mathematical research. We consider all of these to
be successful outcomes. We also had some unsuccessful outcomes where
students stalled at the literature search phase or had outside commitments preventing them from progressing.

\subsection{Embedded Tutors}
Each semester at \RU there are one or two Calculus II
class sections, each with its own embedded tutor. Each section instructor
informally trains their own tutor. Typically one of the instructors is
experienced with this course and can act as a secondary faculty resource.
We intend for tutors to attend all classes, hold regular office hours, test
out the computer labs ahead of time, and work with groups both inside
and outside of class. In practice, we often are unable to find a qualified
student whose schedule allows them to attend all class periods, and so
we loosen the requirement to attending at least once per week. The use
of the tutor varies by instructor. Some have students working examples
in class, and they use the embedded tutors as an extra help. Other
instructors only lecture. Students in sections where the embedded tutors
help during the class period appear to be more likely to work with the
tutors outside of class. Some tutors try the class's computer assignments
ahead of time. But the tutors are not always given assignments prior
to class, as the assignments are sometimes still being written the night
before the class meets. The tutors do help out during class periods when
computers are used. Tutors are not needed as graders, as the homework
is online. Instructors grade weekly quizzes by hand, as this helps them
to have a sense of where the class is mathematically. Instructors also grade the project parts.
Tutors are student workers paid hourly; their salary is part of the
institutional budget, often including federal work study.

\subsection{Project Designers}
At \RU, many students transfer in or take calculus their sophomore year, which means they are ready for a traditional undergraduate research experience only in their senior year. Therefore there is a need for students here to have research experiences requiring less background knowledge.  
Project creation allows student researchers the chance to apply calculus to an area that interests them, rather than just following faculty interests.  This gives them a wide range of choices.  We ask them to look at an applied problem involving actual data, ideally with some sort of social justice component.  In the initial course redesign process, we had research students do a literature search on calculus projects.  This, along with projects that have been used in our calculus classes here, gives our project designers several distinct examples to look at along with a basic structure to work from without being overly prescriptive.  Althought they were mathematically constrained to design a modeling project for a calculus class, designers independently explored and selected the application that would be used. We meet weekly with our research students to talk about their progress, ideas, and challenges.  During the week, they work independently, although we are always available either in person or by e-mail.    

Students are funded through a STEP grant shared with the sciences and through our university's honor's program.

\section{CURRICULAR DESIGN AS STUDENT RESEARCH\label{curricular}}

In the report {\it Mathematics Research by Undergraduates: Costs and
Benefits to Faculty and the Institution} \cite{MAA}, the Committee on the Undergraduate Program in Mathematics of the Mathematical Association of America lists four characteristics of undergraduate mathematics research:
\begin{itemize}
\item The student is engaged in original work in pure or applied mathematics.
\item The student understands and works on a problem of current research interest.
\item The activity simulates publishable mathematical work even if the outcome is
not publishable.
\item The topic addressed is significantly beyond the standard undergraduate curriculum.
\end{itemize}
 
Our research students create projects for use in a Calculus II classroom,
and so theirs is more of an applied curricular design research project
than a traditional mathematics research project. Because of this, the
first item is only partly true; the work is often adapted for a Calculus II
classroom from another source. The second item holds, and was a significant motivator for our research students when they chose the topic
of their project. The third holds in the sense that their work, when
completed, is published and used in our classrooms. For our students,
two of six projects reached this point. Others either lacked time (or
good data sets) or morphed away from a Calculus II project and into applied math research for an honor's thesis. The final point applies in the sense that it takes them outside the traditional curriculum. While the mathematics might be found in an undergraduate math modeling
course, the act of designing math curricula and relating math to a social
justice theme provide a deeper challenge. At the same time, this allows
our student project designers the chance to work on research very early
in their careers here.

A tiered description for the student research work comes from \cite{Olaf}.  In it, they define three levels of research activities:

``{\bf Guided discovery:} In these classroom activities, students make step-by-step
progress toward a standard (but unknown to them) mathematical formula, or
other result, via open-ended, but guided questions.

{\bf Independent investigation:} In these multi-day activities, the instructor asks
open-ended questions that require independent exploration by the students.
Results may not be surprising to professionals, but they cannot be easily found
in the literature.

{\bf Scholarly inquiry:} In these intense activities, students engage in scholarly
work that is typical of a given field of inquiry.''

Our research students engage in scholarly inquiry in curriculum design, researching applied areas and educational theories in order to develop a guided discovery-style project for the Calculus II class itself to
use. To the extent they are trying to develop an original calculus project
they are also engaged in independent investigation. One of the designers
on his own created a template for studying regression. This creates a
cycle of learning, where our more advanced students progress beyond
their Calculus II course while at the same time helping students who are in
the Calculus II course itself. We have found a few instances where advanced
students created mathematics materials for introductory students in the
literature. In \cite{Partners}, four mathematics students worked with lecturers to create materials for a module in vector spaces and complex variables. The authors noted the call for student-led curricular design in the UK \cite{Call1,Call2}, which has been responded to in other fields. The authors also
noted that there was a paucity of literature on student-created mathematics curricula. At least two papers were written in response to \cite{Partners}. In \cite{MorePartners}, two students at an English university worked to create mathematical modeling teaching and assessment tasks for a second-year mathematics for engineers course.  In Swinburne University of Technology in Australia \cite{Videos}, a team of engineering and multimedia students created videos for engineering students to demonstrate how mathematics is used in engineering.  

In \cite{AEW}, the program of Academic Excellence Workshops at Cal Poly
Pomona involved STEM upperclassmen as leaders of cooperative-learning
based workshops for underclassmen in courses ranging from college algebra through calculus. In it, they led weekly problem sessions, and the student facilitators selected the materials used. The facilitators met
weekly with faculty who were teaching the course, and they went to
an intensive two-day training. Although the paper does not mention
whether the problems are student-created or student-selected, the process of choosing appropriate course materials is an advanced one, and so
this is a notable example of students contributing to the enhancement of mathematics curricula.
Some institutions involved with the Science Education for New Civic Engagements and Responsibilities (SENCER) project are also working with students to create curricular materials, notably in biology \cite{biologists} where students are used to create and update labs.  At Guilford College students are creating a new course as a part of their independent study \cite{guilford}, and at New England College a proposal is being piloted \cite{newengland}.  At the United States Military Academy students are doing in-depth assessment research of the university's curricular design across the STEM disciplines \cite{westpoint}.

In many of these cases, a small number of students were selected to
participate in this work, but without a particular common experience
to draw upon. In our project the induction process is systematic and
intentional and provides the cycle of learning described above. Students
have the initial experience of working on a Calculus II project as students
in the class, then are given the opportunity to work as a peer tutor
or project designer (or both). Their subsequent work then impacts the
next set of potential tutors and designers. The depth of the creative
work done by our project designers appears to be beyond that of the
AEW leaders, and so the combination of a cycle of learning with the
depth of experience appears to be unique to our endeavor.

\section{METHODOLOGY\label{methodology}}
This qualitative study involves a relatively small number of participants: five project designers, one of whom also was a tutor, and eight additional students who were embedded tutors.  We proceeded by interviewing each student with several open ended questions (Appendix A) intended to get them to reflect on how they were affected by the experience. We created a survey after we interviewed a few of the students, and it included questions that were based on the interviews.  The survey itself was anonymous, and it was used to corroborate the interviews. Eight
students agreed to be interviewed; four of these also completed a follow
up survey. Two individuals completed the survey, but not an interview.
Three did not respond to our contact request. 

To categorize the responses, the three authors independently reviewed the interview transcripts and labeled responses according to a variety of categories (Appendix B). The labels were compared and discussed until reaching consensus. The results are organized into three main themes as follows:\\
{\bf Theme A:} Insight into better learning processes.\\ 
{\bf Theme B:} Insight into applying mathematics/calculus.\\ 
{\bf Theme C:} Feedback on improving the experience of embedded tutors and researchers. 

\section{EXPERIENCES AND RESULTS\label{results}}
\subsection{Experiences}
Overall, our experiences have been positive. In the first part of this section we will describe some of our observations made as course instructors and research advisors. In the second part of the section we will concentrate on the actual results of our interviews and surveys.

Our embedded tutors have had varied experiences; some only gained a review of calculus, while others developed into expert teachers.  But the students surveyed agreed with us that they all gained in some form, regardless of the amount.

At the beginning, we hoped that the use of tutors would help to foster a sense of community among the students in the class and in our major.  We also hoped that the class's mathematical skill level would increase as well as the tutor's mathematical skills.  We hoped for smoother computer labs, smoother group dynamics during the project, and a source of peer advice. Two of the tutors explicitly commented on the increased sense of community; we observed this as well, both in the classroom and among
the tutors. It was difficult as instructors to discern whether there was
an effect on the mathematical skill level of the class or on their group
dynamics. This is largely due to the small number of class sections
observed. But there was a noticeable effect on the computer labs; these
benefited greatly from the extra support. The amount of peer advice
given varied by tutor; some of them commented on this in the interviews.

There has not been a good mechanism for the class to give feedback on the tutors --an online survey had a low response rate-- but informally they compliment tutors who are actively involved.

Our experiences with student researchers has also been mixed. They
have definitely learned the difficulty of finding data; since much of what is
online is processed data which gives only means, medians, and standard
deviations rather than raw data. They usually found that government
sites are a good data source. As a result of their work, we have used
two student-created projects in our course; these are on modeling population and modeling air pollution. Those student researchers gave talks
on their projects, both internally and externally. We also had student
work develop into more involved research projects that were not used in
class. These were on actuarial models and modeling head injuries. The
students working on these areas gave talks on them internally (actuarial)
and externally (head injury). Two projects (population, actuarial modeling) developed into honors theses, with the first thesis also studying
the impact of the population project on the class using it. We have had
some projects that were not finished. One of the student researchers,
working on temperatures, was stalled in the data collection stage, and
did not come up with a way to relate the topic to calculus. The other,
working on planetary motion, had some activities planned, but they were
lost in a move and not typed up. After this, we started making students
type up their results part-way through their research project. This helps
to prevent the loss of work.

In our experience, project designers have the best results when they fill out weekly timesheets in order to be paid for their work, rather than being paid in a lump sum for their summer research.  We think they feel more accountable and are better able to pace themselves with the timesheets.

\subsection{Results}
The student interviews provide indications that the students benefited from their experience as tutors and designers as well as from working on the Calculus II projects. They also provide valuable feedback on the curricular design. 

Note that we have removed words such as ``Uh, um, like'' as well as repeated phrases from the transcription quotes without explicitly labeling each instance of this.

\subsubsection{Theme A: Insight into better learning processes}
This theme encompasses the students' sense of themselves as learners and as tutors, how math instruction is enhanced by students working on open ended problems, and the components of effective project design. All of the tutors and designers report some gains in their understanding of calculus and in becoming better students themselves. All appreciate the value of having a project required in Calculus II. Tutors and designers put a lot of thought into what students need to be successful.

Tutors noted the value of sitting through the class a second time. They were able to work on areas of weakness and were better able to look for connections among the topics and applications. Having experienced the challenge of working on a project that is more open-ended than a typical homework problem, they are in a position to coach students through the process. One tutor spoke at length about the psychology of a student facing a difficult subject. The self-doubt the student feels was also experienced by the tutor, but knowing that the struggle is shared by the tutor can help to get the student through to the other side. 

Project designers tried to include elements that connected naturally to particular calculus concepts. For example, population growth naturally associates with differential equations. But more importantly they tried to make the project connect to students' own majors such as biology. The project designers discussed how they had to think about what calculus topics students needed to know and how the project can help them with difficult concepts. One project designer explained that conceptually, integration is difficult for students, and so he wanted the project to connect integration to a real life problem. They are interested in making the topics current, not ``like 18th century.'' By putting more emphasis on a meaningful situation, students would naturally move away from a more mechanical view of calculus. 

Several tutors viewed the project as something that motivates interest in math. Previously they had experienced math classes as memorization and refinement of processes. As an embedded tutor they could appreciate the value of the teacher setting a context relevant to the math. \\
One said, ``I think that it was really interesting getting to do lots of different things, but I also think that it is something that students talk about especially within the same degree program. So if we did something that was more biological, population based... one semester when I had a classmate who did something that was more ecological, like the oil spill one, we could have those conversations about how we’re applying the same skills in a very sort of different context.''

It is evident that tutors and creators think a lot about the students. They care about whether the project is feasible and relevant to student interests. The majority of the students in Calculus II are science majors, so project designers looked for projects that related to biology and chemistry.   Several tutors expressed empathy for the students and were motivated to help students practice, find related problems in the homework, and discover new ways to explain things.

Tutors took advantage of the unique relationship they have with the students as an embedded tutor. Tutors know what the students are hearing from the instructor; tutors can fill in gaps from the instructor to the students and can also give some of the students' perspective back to the instructor. This advocacy for the students helps the instructor better understand the needs of the students. The tutor's view is different from the instructor, principally that they have only recently learned the course material and have a better sense of where the students' thinking is presently.  Students often felt more comfortable talking to a peer. 

\subsubsection{Theme B: Insight into applying mathematics/calculus}
Our main motivation for incorporating projects in Calculus II is to give students the ability to talk about calculus and its uses. Our tutors and designers reflected about their time as a calculus student in their interviews.  The project challenges students to think about the mathematical concepts in a contextualized situation that requires imagination and technological assistance.  Calculus II students must communicate among themselves about mathematical modeling in order to successfully complete the project.
Many cited this communication as crucial. One described how they had done group work in their previous calculus class but ``it was never actually going out into the world and presenting your findings and being knowledgeable of what you were talking about, so I liked that as a component.'' One also noted how it helped them talk to professionals at a job fair.

The project designers' reflections deepened when discussing the thinking that went in to designing a project.  Project design required student researchers who were themselves Calculus II students a semester or two ago  to look for ideas that were feasible for Calculus II students to complete in a semester. Designers wanted their projects to be socially relevant,
  
so they searched for an interesting area and then had to deconstruct it: one chose to study head injuries and came across the head injury index. That led to a new kind of analysis for her, working backwards from a formula to work out its derivation.

The designers intended for students to experience how a model may be limited, but they still wanted students to make valid inferences about what formulas would be reasonable to try. One designer noted his own growth as a student through understanding why concepts are true over just accepting them as an established principle.
 
The project designers were able to apply knowledge acquired since having had Calculus II. One, an actuarial science major, designed a project using mortality tables. Reflecting on the project done and the project design led to the problem of data. The projects needed some publicly available data to analyze. They could see that the data used when doing the project as a Calculus II student had problems. Most of the designers expressed awareness of the difficulty of doing a project with real data, in particular, finding a good source and dealing with flaws in the data itself. 
There is consensus among the designers that the project brings value to the class. It gives insight into how calculus can be applied in the real world, and the learning that is needed to navigate the project provides an incentive for students to learn more about calculus itself.

\subsubsection{Theme C: Feedback on improving the experience of embedded tutors and researchers.}

Interviewing the tutors and researchers brought the opportunity for some feedback on how to run the different activities. The tutors felt strongly that more preparation and better coordination between instructors and tutors was needed. They gave some suggestions about the structure of the class, and some insights on the value they should bring to it. 

Tellingly, the project designers did not express concerns about what was expected of them. Their biggest concern regarded the difficulties of finding good projects, particularly those with usable data sets. Because the designers met regularly with their research mentor, they remained informed of the goals and expectations of the project.

Most tutors saw the value and importance of integrating technology into the class, but most did not feel that they improved their skills during their participation as tutors. In fact many pointed out the need for more training for students, tutors and instructors. The tutors believe that students in the class need more formal instruction on using the software, since they feel that much time in class is spent troubleshooting the difficulties students are having or getting them started. The tutors felt that more training for them would help them be more effective at helping the students, as they were unable to answer some questions students had. Finally, there are indications that it is mainly the instructors that need additional training, both on the software being used and on the way to utilize the tutors effectively during the semester. In some cases the instructor relied on the tutor to troubleshoot any problems arising with the software. Most tutors felt instructors only explicitly engaged them when technology was being used in that day's class. In addition, many of the tutors were not active during class unless there was an activity involving computers. 

It is not surprising then that communication was the most cited concern among tutors. Several of them said they wished they knew more about the instructor's goals. The true value of the embedded tutor is to act as a partner of the instructor, and for this he/she needs to be aware of what the instructor is trying to accomplish. Some tutors tended to hold back and not be proactive about helping, in part because they had no direction and in part by their own inexperience and lack of training. Many noted the value of having the time structured so that tutors are available to students both in and outside of class. Opportunities to be active in the class were important to the tutors though some needed more prompting from the instructor.

This suggests some changes in the structure would help facilitate the tutor's activities. Possibilities include more training involving all members of the team, regular meetings between tutor and instructor where plans for the class are discussed, and a set of prompts for the instructor to help guide the tutor.

\section{CONCLUSIONS}

Our experiences with student researchers mirrored those of others, even though our student research had a curricular focus, instead of a mathematical one. In \cite{ThreeYearStudy}, a survey of 76 student science researchers at four different liberal arts institutions was compared with literature from 54 different papers on hypothesized benefits of being a student researcher. They found that students reported gains in many areas, including confidence in their ability to do research, finding connections between and within science, their organizational and computer skills, their enthusiasm, enhanced resumes, and their attitudes towards learning and working as a researcher.  In our study, we also found these, giving evidence that this type of student research project has many of the same benefits of a traditional research project.

We found that curricular design is a valid research experience for undergraduates, leading to many of the benefits of a mathematically focused project. The main advantage of this type of research with a curricular focus instead of a disciplinary one is that it is possible for students to do this type of work when they are just beyond the calculus level. In our study, designers and tutors gained a deeper knowledge of how to apply mathematics and use technology.  Both reflected on what makes a good teacher, indicating this type of experience could greatly benefit undergraduates who are interested in teaching.  They also provided thoughtful comments on how to improve the program.

It would be interesting to have introductory-level students  designing work for earlier courses, such as college algebra and trigonometry.  This might provide a way for students who come into college under-prepared to connect with a mathematics major in a meaningful way.  We would be interested to see work in this direction.

\section*{ACKNOWLEDGEMENTS}
We would like to thank Amy Dexter, Sherri Berkowitz, and Amanda Fisher for giving us pointers on the qualitative research process.  We would like to thank the SENCER project for the initial impetus to redesign our Calculus II course.  Thank you to the NSF STEP grant and the honor's program for supporting the student researchers financially, and to the Provost's office for providing some travel suport and support for a student worker.  Thank you to Janet Campos for her work transcribing the interviews.

\section*{APPENDIX ~A: INTERVIEW QUESTIONS}
\begin{itemize}
\item Describe your experiences as a student in Calculus II. 

\item (Designer only) How did you go about creating the Calculus II project? 
\item (Designer only) What did you learn while creating a Calculus II project?

\item (Tutor only) What did you do as an embedded tutor?
\item (Tutor only) What did you learn as an embedded tutor for Calculus II?

\item What do you think an embedded tutor should do?

\item If you could travel back in time, what advice would you give to yourself?

\item What resources would be useful for you to have?

\end{itemize}

\section*{APPENDIX ~B: CODING CATEGORIES}
\begin{itemize}
\item Teaching style of instructors:  
How it influenced you when you took calculus
Learning in class you were taking 
Teaching style of instructor in the class you were taking
Teaching style of instructor in the class you were tutoring

\item Teaching style of self:
Your learning due to you acting as a tutor 
Student learning due to your tutoring
Communicating with students you were tutoring
	Project creation: reflecting about how to help other students

\item Resources:
Communicating with faculty
Resources (computer-based)
Resources (other)
``People creating/tutoring should have property X.''
``I needed knowledge about X.''
``I used X to do Y.''

\item Uses of calculus:
Uses of calculus (class you were taking)
Uses of calculus (class you were tutoring)
Uses of calculus (project you created)

\item Self-reported changes:
In how you think about calculus
``I created a project, and now I am awesome at calculus''
Interest in mathematics/applications/teaching
Communicating with an outside audience
``I can do X.''
\end{itemize}

\section*{APPENDIX ~C: SURVEY QUESTIONS}
\begin{itemize}
\item When did you master calculus? (Select all that apply.) \begin{itemize}
 \item As a student in calculus class.
 \item As a tutor.
 \item While creating a new project (as my research).
 \item In Calculus III.
 \item In more advanced math or actuarial science courses.
 \item In my job.
 \item Other: 
\end{itemize}
\item What background should a student have before working on designing a new project? (Select all that apply.)\begin{itemize}
 \item Calculus II.
 \item Calculus III.
 \item More advanced mathematics courses.
 \item A mathematical computing course.
 \item A programming course.
 \item Other:  
\end{itemize}
\item What software have you used (select all that apply from Maple or Mathematica, Wolfram Alpha, Excel, Powerpoint, Statistical software (specify type), Other) \begin{itemize}
\item at \RU outside of Calculus II?
\item at \RU while taking or tutoring Calculus II? 
\item at your job, if you are employed outside of \RU?
\item What software do you think is useful for Calculus II students to learn?
\end{itemize}
\item Ignoring exam days, about how frequently should a Calculus II class involve (select from Never, Once or twice, Monthly, Weekly, Every day for each option below.)\begin{itemize}
\item A lecture.
\item Group work.
\item Project work.
\item Computer work.
\item Problem sessions.
\end{itemize}
\item How do projects benefit Calculus II students? (Select from Essential for the students, Helpful for the students, Somewhat helpful for the students, They'd learn this without the projects for each option below.)\begin{itemize}
\item Increased computational skills (taking integrals, etc).
\item Increased conceptual skills (what does an integral represent, etc.).
\item Finding out that calculus is actually useful.
\item Increased skills working with data.
\item Increased skills working with people.
\item Increased skills with computers.
\item Learning about connections to other fields.
\item Better able to communicate mathematics.
\end{itemize}
\item If you created a project, how did constructing a project benefit you? (Select from Essential for the students, Helpful for me, Somewhat helpful for me, Did not help me/does not apply for each option below.)\begin{itemize}
\item Increased computational skills (taking integrals, etc).
\item Increased conceptual skills (what does an integral represent, etc.).
\item Finding out that calculus is actually useful.
\item Working with data.
\item Working with people.
\item Working with computers.
\item Learning about a field you are interested in.
\item Increased patience.
\item Understanding how people think and learn.
\item Understanding how you think and learn.
\item Better research habits.
\item More experience doing literature searches.
\item More responsibility.
\end{itemize}
\item What skills should an ideal embedded tutor have? (Select from Essential, Good to have, Not needed for each option below.)\begin{itemize}
\item Strong mathematical content knowledge.
\item Knowledge of computer software.
\item Love of mathematics.
\item Desire to teach.
\item Enthusiasm.
\end{itemize}
\item What should an ideal embedded tutor do? (Select from Essential, Good to have, Not needed for each option below.)\begin{itemize}
\item Act as a bridge between faculty and students.
\item Help with project during class.
\item Help with calculus examples during class.
\item Help outside of class.
\item Make connections between the project and the class.
\item Help build student confidence.
\item Meet regularly with faculty outside of class.
\end{itemize}
\item Which skills increased for you due to your work as an embedded tutor or project designer? (Select all that apply.) \begin{itemize}
 \item Presenting posters.
 \item Giving talks.
 \item Writing papers/projects.
 \item Tutoring.
 \item Confidence in my ability to do math.
 \item Confidence in my ability to explain math.
 \item Other: 
\end{itemize}
\item What resources are needed for embedded tutors? (Select from Essential to have, Helpful to have, Somewhat helpful to have, Not needed/not applicable for each option below.) \begin{itemize}
\item MyMathLab Course Access.
\item Tutoring time in 401/Math and Science Resource Center.
\item A room dedicated to only tutoring (not classes).
\item Regular communication with faculty teaching the course.
\item Tutor orientation and training prior to the start of the semester.
\item Computer orientation and training.
\item Written sample solutions for homework assignments.
\item Written sample solutions to project parts.
\item Ability to schedule office hours each week at times that are convenient for both you and students (not necessarily fixed).
\item Feedback on performance.
\end{itemize}
\item If you would like to clarify or expand on any of your answers, please do so below.
\end{itemize}

\section*{BIOGRAPHICAL SKETCHES}

S. Cohen is an associate professor of mathematics at \RU.  
He teaches courses to both majors and non-majors throughout the curriculum with particular interest in the History of Mathematics and Abstract Algebra. He is a member of the steering committee of the Chicago Symposium on Excellence in Teaching Undergraduate Mathematics and Science. He earned an M.S. and Ph.D. in Mathematics from University of Illinois Chicago and served as a visiting assistant professor at Loyola University of Chicago.  Steve likes to play undisclosed games of uncertain outcomes.  He also bakes an excellent cheesecake whose outcome is much more certain.

\vspace*{.3 true cm} \noindent B. Gonz\'alez-Ar\'evalo is an associate professor of mathematics at Hofstra University. Previously she was an associate professor of mathematics, statistics, and actuarial science at \RU. Her current research interests include Statistics, Applied Probability and the Scholarship of Teaching and Learning Mathematics. She earned an M.S. and Ph.D. in statistics from Cornell University, and worked as an assistant professor at the University of Louisiana at Lafayette.  She enjoys baking and has two beautiful boys.  It is important to note that she does not bake the boys.

\vspace*{.3 true cm} \noindent M. Pivarski is an associate professor of mathematics at \RU. She is currently serving as the department chair for mathematics and actuarial science. Her current research interests involve heat kernels and their applications in metric measure spaces. Recently, she has been inspired to include students in her research work. This led her to work in the scholarship of teaching and learning mathematics. She earned a Ph.D. in mathematics from Cornell University and worked as a visiting professor at Texas A\&M University.  She likes to eat her co-author's creations, as she is too busy chasing her toddler to bake on her own.

\end{document}